# A test for equality of multinomial distributions vs increasing convex order

## Arthur Cohen[1],[*] John Kolassa[1],[†] and Harold Sackrowitz[2],[*]

*Rutgers University*

**Abstract:** Recently Liu and Wang derived the likelihood ratio test (LRT) statistic and its asymptotic distribution for testing equality of two multinomial distributions vs. the alternative that the second distribution is larger in terms of increasing convex order (ICX). ICX is less restrictive than stochastic order and is a notion that has found applications in insurance and actuarial science. In this paper we propose a new test for ICX. The new test has several advantages over the LRT and over any test procedure that depends on asymptotic theory for implementation. The advantages include the following:

(i) The test is exact (non-asymptotic).

(ii) The test is performed by conditioning on marginal column totals (and row totals in a full multinomial model for a $2 \times C$ table).

(iii) The test has desirable monotonicity properties. That is, the test is monotone in all practical directions (to be formally defined).

(iv) The test can be carried out computationally with the aid of a computer program.

(v) The test has good power properties among a wide variety of possible alternatives.

(vi) The test is admissible.

The basis of the new test is the directed chi-square methodology developed by Cohen, Madigan, and Sackrowitz.

## 1. Introduction

Recently, Liu and Wang [8] derived the likelihood ratio test (LRT) statistic and its asymptotic distribution for testing equality of two multinomial distributions vs. the alternative that the second distribution is larger in terms of increasing convex order (ICX). See also Liu and Wang [7]. A formal definition of ICX is as follows: the distribution of a random variable $Y$ is larger than the distribution of a random variable $X$ in the increasing convex order, i.e. $X \leq_{ICX} Y$, if and only if $E\{f(X)\} \leq E\{f(Y)\}$ holds for all non-decreasing convex functions $f$ for which expectations are defined. ICX is less restrictive than stochastic order and is a notion that has found applications in insurance and actuarial science. See, for example, Goovaerts, Kaas, Van Heerwaarden and Bauwelinckx [6] and other references cited by Liu and Wang [8]. In this paper we propose a new test for ICX. The new test has several advantages over the LRT and over any test procedure that depends on asymptotic theory for implementation. The advantages include the following:

[1]Department of Statistics, Rutgers University, 110 Frelinghuysen Road, Piscataway, NJ 08854-8019, e-mail: artcohen@rci.rutgers.edu; kolassa@stat.rutgers.edu

[2]Department of Statistics, Rutgers University, 110 Frelinghuysen Road, Piscataway, NJ 08854-8019, e-mail: sackrowi@rci.rutgers.edu, url: www.stat.rutgers.edu

[*]Research supported by NSF Grant DMS-0457248.

[†]Research supported by NSF Grant DMS-0505499.

(i) The test is exact (non-asymptotic). It can be implemented regardless of sample sizes.

(ii) The test is performed by conditioning on marginal column totals (and row totals in a full multinomial model for a $2 \times C$ table). Conditioning enables the null hypothesis to be expressed as a simple null, and can be carried out by calculating conditional $P$-values.

(iii) The test has desirable monotonicity properties. That is, the test is monotone in all practical directions (to be formally defined). Intuitively monotone in practical directions means that if the test rejects for a sample point, say $\mathbf{x}$, then it should also reject for a sample point $\mathbf{y}$ where $\mathbf{y}$ empirically is more indicative of ICX than $\mathbf{x}$. The LRT is not monotone in all practical directions.

(iv) The test can be carried out with the aid of a computer program.

(v) The test has good power properties among a wide variety of possible alternatives.

(vi) The test is admissible.

The basis of the new test is the directed chi-square methodology developed by Cohen, Madigan and Sackrowitz [5].

In the next section we will state the formal model while defining ICX. We will also state the hypothesis and define practical directions. Furthermore, we determine the practical directions for the ICX alternative. In Section 3, we offer the directed chi-square test statistic. Section 4 contains an example concerned with an age discrimination study. In this same section we offer a simulation study comparing powers of the new test with an exact version which uses the LRT statistic. Finally Section 5 contains a discussion regarding the importance of the monotonicity properties.

## 2. Models and definitions

Consider a $2 \times C$ contingency table under the product multinomial model. Assume the $C$ categories are ordered (worst to best; increasing age groups; etc.). Let $X_{ij}$, $p_{ij}$, $i = 1, 2$; $j = 1, \ldots, C$ represent cell frequencies and cell probabilities for the cell $i$, $j$. Note $\sum_{j=1}^{C} X_{ij} = n_i$ are fixed, $\sum_{j=1}^{C} p_{ij} = 1$ for $i = 1, 2$, and let $X_{1j} + X_{2j} = t_j$ denote column totals, $j = 1, \ldots, C$. Also let $N = n_1 + n_2$. Define log odds ratios as

$$(2.1) \qquad \nu_j = \log(p_{1j} p_{2C} / p_{1C} p_{2j}), \quad (j = 1, \ldots, C-1).$$

Also $\mathbf{X} = (\mathbf{X}_1, \mathbf{X}_2)'$, where $\mathbf{X}_i = (X_{i1}, \ldots, X_{iC})$, for $i = 1, 2$, $\nu = (\nu_1, \ldots, \nu_{C-1})'$. Note $\mathbf{X}$ is a $2C \times 1$ column vector and $\nu$ is a $(C-1) \times 1$ column vector. The null hypothesis to be studied is $H : \mathbf{p}_1 = \mathbf{p}_2$. The alternative hypothesis is called increasing convex order (ICX) and is defined as follows: Let $\lambda_1 > \lambda_2 > \cdots > \lambda_{C-1} > 0$ be $(C-1)$ given constants. Then the distribution with parameter $\mathbf{p}_1$ is said to be smaller in ICX than the distribution with $\mathbf{p}_2$ as parameter if for $r = 1, \ldots, C-1$,

$$(2.2) \qquad \Delta_r = \lambda_r \sum_{j=1}^{r} (p_{1j} - p_{2j}) + \sum_{j=r+1}^{C-1} \lambda_j (p_{1j} - p_{2j}) \geq 0.$$

This definition is essentially the same as the one given by Liu and Wang [8]. This is an equivalent form of ICX for two multinomial distributions. Hence the alternative is denoted by

$$(2.3) \qquad K_{ICX} : \big\{ \mathbf{p} = \mathbf{p}_1, \mathbf{p}_2 \big) : (2.2) \text{ holds} \big\} \setminus H.$$

Our approach to testing is to condition on the column totals (row totals as well if the model is full multinomial) since these totals are the complete sufficient statistics under $H$. We let $\mathbf{m} = (t_1, \ldots, t_C)$ denote these sufficient statistics (under the full multinomial model $\mathbf{m} = (n_1, n_2, t_1, \ldots, t_C)$). The conditional distribution of $\mathbf{X}^{(1)} = (X_{11}, \ldots, X_{1(C-1)})'$ given $\mathbf{m}$ is the multivariate extended hypergeometric distribution, which in exponential family form is

$$(2.4) \qquad f(\mathbf{x}^{(1)}; \nu) = h_{\mathbf{m}}(\mathbf{x}^{(1)}) \beta_{\mathbf{m}}(\nu) e^{\mathbf{x}^{(1)'} \nu}.$$

See [3]. For the conditional problem, $H$ becomes $H^* : \nu_1 = \nu_2 = \cdots = \nu_{C-1} = 0$. In order to specify the appropriate alternative when $\mathbf{m}$ is fixed we need

**Lemma 2.1.** *Let $Q^- = \{\nu \in \mathbb{R}_{C-1} : \nu_j < 0, \text{ all } j = 1, \ldots, C-1\}$. Consider the set*

$$(2.5) \qquad \Gamma = \{\nu \in \mathbb{R}_{C-1} \setminus Q^-\} \setminus \{\mathbf{0}\}.$$

*Given any $\nu \in \Gamma$, there exists some $\mathbf{p}(\nu)$ satisfying (2.2). Furthermore if $\nu \in Q^-$, there is no $\mathbf{p}$ satisfying (2.2).*

*Proof.* See Appendix. □

In light of Lemma 2.1, for the conditional problem we take the alternative to be

$$K^*_{ICX} : \{\nu : \nu \in \Gamma\} \setminus H^*.$$

Now let $\phi(\mathbf{x})$ denote a test function; i.e., $\phi(\mathbf{x})$ is the probability that the test rejects $H$ for an observed sample point $\mathbf{x}$.

**Definition 2.1.** A test $\phi(\mathbf{x})$ is said to be monotone in the direction $\xi = (\xi_1, \ldots, \xi_{2C})'$ if and only if

$$(2.6) \qquad \phi(\mathbf{x}) \leq \phi(\mathbf{x} + \gamma \xi),$$

for every $\gamma \geq 0$.

Since we will do testing by conditioning on $\mathbf{m}$, and since $n_1, n_2$ are fixed, hereafter we only consider directions such that

$$(2.7) \qquad \xi_{1j} + \xi_{2j} = 0, \quad j = 1, \ldots, C, \quad \text{and} \quad \sum_{j=1}^{C} \xi_{ij} = 0, \quad i = 1, 2.$$

At this point let $\hat{p}_{ij} = x_{ij}/n_i$ and consider the vector $\hat{\mathbf{p}} = (\hat{p}_{11}, \ldots, \hat{p}_{2C})'$. Let

$$(2.8) \qquad \Delta^*_r(\mathbf{x}) = \Delta_r(\hat{\mathbf{p}}).$$

**Definition 2.2.** A direction $\mathbf{d}$ is said to be a practical direction if

$$(2.9) \qquad \Delta^*_r(\mathbf{x} + \mathbf{d}) \geq \Delta^*_r(\mathbf{x}), \quad \text{for } r = 1, \ldots, C-1.$$

An interpretation of a practical direction is that the empirical distributions are becoming more ICX. Note that if a test function is monotone in directions $\mathbf{d}_1$ and $\mathbf{d}_2$ (see (2.6)), it is monotone in the direction $a_1 \mathbf{d}_1 + a_2 \mathbf{d}_2$ as long as $a_1 \geq 0$, $a_2 \geq 0$. This implies that the collection of practical directions for which $\phi$ is to be monotone generates a closed convex polyhedral cone $\mathcal{C}$. Using (2.2), (2.7), (2.8), and (2.9) we may express $\mathcal{C}$ as follows:

$$(2.10) \qquad \mathcal{C} = \{\mathbf{d} : B\mathbf{d} = \mathbf{0}, G\mathbf{d} \geq \mathbf{0}\},$$

where $B$ is a $(C+1) \times 2C$ matrix expressing the constraints in (2.7) and $G = (G_1, -(n_1/n_2)G_1)$ is the $(C-1) \times 2C$ matrix and

$$(2.11) \qquad G_1 = \frac{1}{n_1} \begin{pmatrix} \lambda_1 & \lambda_2 & \cdots & \cdots & \cdots & \lambda_{C-1} & 0 \\ \lambda_2 & \lambda_2 & \lambda_3 & \cdots & \cdots & \lambda_{C-1} & 0 \\ \lambda_3 & \lambda_3 & \lambda_3 & \lambda_4 & \cdots & \lambda_{C-1} & 0 \\ \vdots & & & & & & \vdots \\ \lambda_{C-1} & \cdots & \cdots & \cdots & \cdots & \lambda_{C-1} & 0 \end{pmatrix}$$

**Remark.** The same example used in [2] can be used to demonstrate that the LRT is not monotone in all practical directions.

## 3. Directed chi-square

The directed chi-square statistic was introduced in [5]. The statistic is

$$(3.1) \qquad \chi_D^2(\mathbf{x}) = \sum_{i=1}^{2} \sum_{j=1}^{C} x_{ij}^{*2}/n_i t_j = \inf_{\mathbf{u} \in A(\mathbf{x})} \sum_{i=1}^{2} \sum_{j=1}^{C} u_{ij}^2/n_i t_j,$$

where $\mathbf{u} = (\mathbf{u}_1, \mathbf{u}_2)'$ is a $2C \times 1$ vector, $\mathbf{x}^*$ is the minimizer of the sum on the right-hand side of (3.1) and $A(\mathbf{x})$ is a set in $\mathbb{R}_{2C}$, depending on the data $\mathbf{x}$ and determined by a set of linear equalities and linear inequalities. Namely,

$$(3.2) \qquad A(\mathbf{x}) = \{\mathbf{u} \in \mathbb{R}_{\mathbf{2C}} : \mathbf{B}(\mathbf{u} - \mathbf{x}) = \mathbf{0},\ \mathbf{G}(\mathbf{u} - \mathbf{x}) \geq \mathbf{0}\},$$

where $B$ and $G$ are specified in (2.10), (2.11).

The statistic $\chi_D^2$ can be determined by using an IMSL subroutine called DQPROG. That is, given an observed value of $\mathbf{x}$, call it $\mathbf{x}_0$, determine $\mathbf{x}^*$ of (3.1). Next use the exact method of Pagano and Halvorsen [9] to generate all tables consistent with the given $\mathbf{m}$ and the conditional probabilities under $H$ of these tables. Sum the probabilities of the sample points for which $\chi_D^2(\mathbf{x}) > \chi_D^2(\mathbf{x}_0)$ plus the probabilities of the sample points for which $\chi_D^2(\mathbf{x}) = \chi_D^2(\mathbf{x}_0)$. The total probability is the conditional $P$-value. If this $P$-value $\leq \alpha$, reject $H$.

The directed chi-square test is monotone in all practical directions. A proof of this is given in [5]. The directed chi-square test is admissible. To show this, recognize first that the test for the ICX alternative is admissible for the stochastic order (SO) alternative, which is a smaller parameter set than the ICX alternative. The admissibility for the SO alternative follows from Theorem 4.3 of [3], using the facts that (i) the test is monotone in $x_{11}$ while $\sum_{j=1}^{k} x_{1j}$, $k = 2, \ldots, C-1$, is fixed and (ii) that the acceptance region of the test is convex. See [5].

**Remark.** Should there be several sample points yielding the same value of $\chi_D^2$ it may be helpful to use a backup statistic as is done in [5].

TABLE 1
*Success and age in competition*

| Age     | 20–30 | 31–40 | 41–50 | 51–60 | Totals |
|---------|-------|-------|-------|-------|--------|
| Success | 1     | 6     | 19    | 4     | 30     |
| Failure | 0     | 4     | 11    | 8     | 23     |

## 4. Example and power comparison

Barry and Boland [1] study the relationship of age and successful employment in Ireland. Table 1 contains relevant data.

This is a reasonable example to consider ICX as opposed to stochastic ordering, since apriori one might suspect that older people will have a smaller chance of gaining employment than younger people, whereas at the young age groups you would not expect much of a difference. Using $\lambda_1 = 3$, $\lambda_2 = 2$, $\lambda_3 = 1$ (see (2.2)), we find the conditional $P$-value for this data set using the directed $\chi^2$-test is 0.10539. For the LRT the corresponding $P$-value is 0.16575.

A study was conducted to compare exact conditional power of the directed $\chi^2$-test with an exact test based on the LRT statistic. The study was based on the data from the marginal totals of Table 1 save that the first and second columns were combined. Hence the problem is in terms of a $2 \times 3$ table with marginal column totals of (11, 30, 12) and row totals (30, 23). Calculations were performed using Fortran 90 and the IMSL mathematical library for nonlinear function minimization. We took $\lambda_1 = 2$, $\lambda_2 = 1$. In order to calculate the constrained maximum likelihood estimate under ICX order, IMSL subroutine DL2ONG was used to minimize the likelihood subject to the linear constraints (2.2) and $\sum_{j=1}^{C} p_{ij} = 1$, $i = 1, 2$. The likelihood is simply the product of the two multinomial distributions. In addition to the likelihood the derivative was also provided by DL2ONG in a separate subroutine. The chi-square statistic (3.1) minimized under ICX order was obtained using the IMSL routine DQPROG for minimizing a quadratic form under linear constraints; the constraints are given in (3.2).

In order to calculate the $P$-value of the directed chi-square test and to calculate powers all tables with the same marginal totals as the observed table are enumerated. These tables, and their probabilities conditional on row and column sums, were calculated using the algorithm of Pagano and Halvorsen [9]. A conditional $P$-value was calculated as the sum of table probabilities for which the test statistic was as large as or larger than that observed. Powers were calculated by reweighting the tables using the ratio of the likelihood under the alternative hypothesis to the likelihood under the null hypothesis, and summing the probabilities associated with tables whose test statistics were as large as or larger than those observed. Also the powers were adjusted so that test sizes are exactly 0.05.

Table 2 contains exact powers of the direct chi-square test and the exact test performed conditionally using the unconditional LRT statistic. Various ICX alternatives are considered. We note that the powers of the two tests are comparable. The LRT is slightly better for some alternatives that are further from a null case while $\chi^2$ is preferred for alternatives closer to a null case.

## 5. Discussion

One referee has misgivings about this paper because of our claim that monotonicity in practical directions is an intuitively desirable property. The referee refers to Perlman and Chaudhuri [11] where it is argued that such a property is not compelling and since the likelihood ratio test does not have the property it is undesirable. Our reaction to this has been discussed in some detail in Cohen and Sackrowitz [4], a paper that appears in the same year of the same journal as the paper by Perlman and Chaudhuri [11].

TABLE 2
*Exact powers for the directed $\chi^2$ test and for the LRT alternatives*

| Alternatives | | | | | | Powers | |
|---|---|---|---|---|---|---|---|
| $p_{11}$ | $p_{12}$ | $p_{13}$ | $p_{21}$ | $p_{22}$ | $p_{23}$ | $\chi^2$ | LRT |
| 0.10 | 0.60 | 0.30 | 0.20 | 0.10 | 0.70 | 0.9747 | 0.9940 |
| 0.10 | 0.60 | 0.30 | 0.20 | 0.20 | 0.60 | 0.8941 | 0.8918 |
| 0.10 | 0.80 | 0.10 | 0.20 | 0.30 | 0.50 | 0.9591 | 0.9743 |
| 0.10 | 0.50 | 0.40 | 0.30 | 0.10 | 0.60 | 0.9790 | 0.9823 |
| 0.30 | 0.50 | 0.20 | 0.40 | 0.10 | 0.50 | 0.9640 | 0.9782 |
| 0.20 | 0.30 | 0.50 | 0.30 | 0.10 | 0.60 | 0.7050 | 0.7009 |
| 0.10 | 0.40 | 0.50 | 0.20 | 0.10 | 0.70 | 0.9163 | 0.9161 |
| 0.10 | 0.60 | 0.30 | 0.30 | 0.10 | 0.60 | 0.9863 | 0.9961 |
| 0.50 | 0.40 | 0.10 | 0.60 | 0.10 | 0.30 | 0.9526 | 0.9627 |
| 0.10 | 0.60 | 0.30 | 0.15 | 0.45 | 0.40 | 0.2278 | 0.2205 |
| 0.20 | 0.30 | 0.50 | 0.25 | 0.20 | 0.55 | 0.1934 | 0.1863 |
| 0.20 | 0.40 | 0.40 | 0.25 | 0.20 | 0.55 | 0.4464 | 0.4369 |
| 0.10 | 0.40 | 0.50 | 0.12 | 0.35 | 0.53 | 0.0862 | 0.0813 |
| 0.30 | 0.30 | 0.40 | 0.32 | 0.27 | 0.41 | 0.0701 | 0.0660 |
| 0.10 | 0.40 | 0.50 | 0.15 | 0.30 | 0.55 | 0.1689 | 0.1647 |

Our take on the controversy is as follows: Likelihood inference is the default methodology in much of statistical inference where it is feasible. It has large sample optimality properties that are unsurpassed under very mild conditions. It generally has intuitive appeal as well. However in some order restricted inference problems likelihood inference has competitors that can have intuitive properties that likelihood procedures do not share. We offer one example here, borrowed from Cohen and Sackrowitz [2] and leave it to the reader to judge the intuitiveness of the monotonicity property we claim is desirable. See also a recent paper by Peddada, Dunson and Tan [10] which offers competitors to maximum likelihood estimators.

**Example.** Consider a $2 \times 3$ contingency table under the product multinomial model. Let $X_{ij}$, $i = 1, 2; j = 1, 2, 3$ be cell frequencies and $p_{ij}$ be corresponding cell probabilities. Test $\mathbf{p}_1 = \mathbf{p}_2$ (when $\mathbf{p}_i = (p_{i1}, p_{i2}, p_{i3})$) ($\sum_{j=1}^{3} p_{ij} = 1$), vs $H_1 : \{\mathbf{p}_2 >_{st} \mathbf{p}_1\} \backslash H_0$, where $>_{st}$ means the $\mathbf{p}_2$ distribution is stochastically larger than $\mathbf{p}_1$ i.e., $p_{11} \geq p_{21}$ and $p_{11} + p_{12} \geq p_{21} + p_{22}$ with at least one strict inequality. Note $\mathbf{p}_2 >_{st} \mathbf{p}_1$ implies $\mathbf{p}_2 >_{ICX} \mathbf{p}_1$. Now consider the following two sample points:

Our intuition suggests that the conditional *p*-value (given marginal totals fixed) should be smaller for sample point 1 than for sample point 2. Yet the *p*-value using the likelihood ratio statistic is 0.169 for sample point 1 and 0.019 for sample point 2.

| Sample Point 1 | | | | |
|---|---|---|---|---|
| Group | Worse | Same | Better | Total |
| Control | 5 | 11 | 1 | 17 |
| Treat. | 3 | 8 | 4 | 15 |
| Total | 8 | 19 | 5 | 32 |

| Sample Point 2 | | | | |
|---|---|---|---|---|
| Group | Worse | Same | Better | Total |
| Control | 0 | 16 | 1 | 17 |
| Treat. | 8 | 3 | 4 | 15 |
| Total | 8 | 19 | 5 | 32 |

We feel that blanket statements that claim monotonicities are always desirable or always undesirable should not be made. Considerations of such should be made on a case by case basis.

## Appendix A: Appendix section

**Proof of Lemma 2.1.** Recognize that $\Gamma$ is the set of $\nu$'s such that at least one component of $\nu$ is greater than zero. Let $\nu_q > 0$ for some $1 \leq q \leq C-1$. Now let

$$p_{1j} = M_1 e^{a_j}, \; j = 1, \ldots, C \text{ but not } j = q \text{ or } j = C,$$
$$p_{1q} = M_1 \Delta e^{a_q}, \; p_{1C} = \Delta e^{a_C},$$
$$p_{2j} = M_2 e^{b_j}, \; j = 1, \ldots, C \text{ but not } j = q \text{ or } j = C,$$
$$p_{2q} = M_2 \Delta e^{b_q}, \; p_{2C} = \Delta e^{b_C}.$$

The constants $\mathbf{a} = (\mathbf{a_1}, \ldots, \mathbf{a_C})$, $\mathbf{b} = (\mathbf{b_1}, \ldots, \mathbf{b_C})$, are as follows:

$$a_j = 3\nu_j/2, \; j = 1, \ldots, q-1, q+1, \ldots, C-1,$$
$$a_q = 0, \; a_C = \nu_1/2,$$
$$b_1 = 0, \; b_j = -(\nu_j + \nu_1)/2, \; j = 2, \ldots, q-1, q+1, \ldots, C-1,$$
$$b_q = -\nu_j - \nu_1/3, \; b_C = 0.$$

This choice of constants yields the given $\nu$'s. The constants $M_1$ and $M_2$ are determined by the fact that $\sum p_{ij} = 1$, $i = 1, 2$. We now verify that this choice of $p(\nu)$ satisfies (2.2) for some $\Delta$. First let $r = 1$, so that we must show

(A.1)
$$\frac{\sum_{j=1}^{q-1} \lambda_j e^{a_j} + \sum_{j=q+1}^{C-1} \lambda_j e^{a_j} + \lambda_q \Delta e^{a_q}}{\sum_{j=1}^{q-1} e^{a_j} + \sum_{j=q+1}^{C-1} e^{a_j} + \Delta(e^{a_q} + e^{a_C})}$$
$$\geq \frac{\sum_{j=1}^{q-1} \lambda_j e^{b_j} + \sum_{j=q+1}^{C-1} \lambda_j e^{b_j} + \lambda_q \Delta e^{b_q}}{\sum_{j=1}^{q-1} e^{b_j} + \sum_{j=q+1}^{C-1} e^{b_j} + \Delta(e^{b_q} + e^{b_C})}.$$

We will let $\Delta \to \infty$ so that from (A.1) it suffices to show

(A.2) $$e^{a_q}(e^{b_q} + e^{b_C}) > e^{b_q}(e^{a_q} + e^{a_C})$$

which reduces to

(A.3) $$e^{a_q + b_C} > e^{a_C + b_q}$$

or

(A.4) $$a_q + b_C > a_C + b_q.$$

However $\nu_q = a_q + b_C - a_C - b_q > 0$ by hypothesis. This shows (3.2) for $r = 1$. For $2 \leq r \leq C-1$ the argument is essentially the same.

To complete the lemma we need to show that if all $\nu$'s are negative then no $\mathbf{p}(\nu) \in \Gamma$. But for $r = C-1$, (2.2) reduces to $p_{2C} \geq p_{1C}$. If this is the case then for some $j$, $j = 1, \ldots, C-1$, $p_{1j} \geq p_{2j}$ implying that $\nu_j \geq 0$.